\documentclass[twoside]{article}
\usepackage{amssymb}
\usepackage{theorem}

\newcommand{\qed}{{\hfill \hbox{\enspace${\square}$}}\smallskip}
\newcommand{\E}{{\mathcal E}}

\newcommand{\R}{{\mathbb R}}
\newcommand{\f}{\varphi}

\newcommand{\Ol}{{\cal O}}

\newcommand{\opn}{{\Ol_{\mathbb P^n}}}

\newcommand{\opt}{{\Ol_{\mathbb P^3}}}

\newcommand{\pn}{\mathbb P^n}
\newcommand{\pnm}{\mathbb P^{n-1}}
\newcommand{\pu}{\mathbb P^1}
\newcommand{\pt}{\mathbb P^3}
\newcommand{\pd}{\mathbb P^2}

\newcommand{\proj}{\mathbb P}

\newcommand{\loc}{\textrm{Locus}}

\newcommand{\om}{\textrm{Hom}}

\newcommand{\pic}{\textrm{Pic}}
\newcommand{\ch}{\textrm{Char}}

\makeatletter
\@addtoreset{equation}{subsection}
\@addtoreset{equation}{section}
\makeatother
\newtheorem{theorem}[equation]{Theorem}
\newtheorem{proposition}[equation]{Proposition}

\newtheorem{lemma}[equation]{Lemma}
{\theorembodyfont{\rmfamily}

\newtheorem{remark}[equation]{Remark}

}

\renewcommand{\theequation}{\arabic{section}.\ifnum\value{subsection}=0\else\ara
bic{subsection}.\fi\arabic{equation}}

\title{\bf{Special rays in the Mori cone \\ of a projective variety}}
\author{Marco Andreatta and Gianluca Occhetta}
\date{} 
\begin{document}
\maketitle

\begin{abstract}
Let $X$ be a smooth $n$-dimensional projective variety over an algebraically
closed field $k$ such that $K_X$ is not nef.
We give a characterization of non nef extremal rays of $X$ of maximal
length (i.e of
length $n-1$); in the case of $\ch(k) = 0$ we also characterize non nef
rays of length
$n-2$.

\end{abstract}
\begin{figure}[bp]
\footnotesize{{\it Mathematics Subject Classification}: Primary 14E30;
Secondary 14J40. }
\end{figure}
\section{Introduction}

Let $X$ be a smooth $n$-dimensional projective variety
over an algebraically closed field $k$ of arbitrary characteristic.

We assume that $K_X$ is not nef, in particular that there exists
an {\it extremal ray} $R$ in the cone $\overline{NE(X)}_{{K_X} <0}$.

The {\it length of the ray $R$} is the integer defined as
$l(R)=$ min $\{ -{K_X}. C: [C] \in R\}$; the set $\loc(R)$ is the set of
closed point $x \in X$ such that there is a curve $x \in C' \subset X$
with $[C'] \in R$.\\
A ray $R$ is said to be {\it nef} if $D. R \geq 0$ for all effective divisors
$D \subset X$, or equivalently if $\loc(R) = X$; if a nef ray exists the
variety $X$ is
therefore uniruled.

\smallskip
The main result of the paper is the following characterization
of the blow-up of a smooth point in terms of extremal rays.

\begin{theorem}\label{main} Let $X$ be a smooth $n$-dimensional projective
variety.
There exists a non nef extremal ray $R$ of length
$\geq (n -1)$ if and only if there exists a morphism
$\f : X \to X'$ into a smooth projective variety $X'$ which is the
blow-up of $X'$ at a point.
\end{theorem}

Let us point out that in the case $\ch(k) > 0$ the existence of a map
$\f: X \to X'$ which contracts all the curves in an extremal ray has not been
proved in general; if $\ch(k)= 0$ this is the so called Kawamata-Shokurov
contraction theorem.

\smallskip
A straightforward application of this result is in the so called {\it
adjunction theory}:
namely if $L$ is an ample line bundle on a a smooth $n$-dimensional
projective variety
$X$ one can study the ampleness of the adjoint bundle $K_X+ (n-1)L$:
this is ample unless there exist extremal rays of length
$\geq (n-1)$.\\
If these rays are nef they will determine the uniruled
variety $X$;
this is the case if $l(R) > (n-1)$, see \cite{Kako}, or if $\ch(k)=0$
is zero, see \cite{Fu1}; if $l(R) = n-1$ the precise characterization of $X$
is still open in positive characteristic (see remark (\ref{nef})).

If there are non nef rays one can apply the above theorem and the following
happens:

\begin{theorem}  Let $X$ be a smooth $n$-dimensional projective variety
with $n \ge 3$
and $L$ an ample line bundle on $X$. \\
If there is a non nef ray in the cone
$\overline{NE(X)}_{(K_X+(n-1)L) \leq 0}$ then,
with the only exception given by the blow-up of $\pt$
in one point, $\pi : Bl_x \pt \to\pt$, and $L = \pi^*(\opt (2))
-[\pi^{-1}(x)]$,
all rays in this cone are non nef and they
can be simultaneously contracted into a smooth variety $X'$, with a morphism
$\f :X \to X'$ expressing $X$ as blow-up of $X'$ at a finite set of points
$B$.\\
Moreover there is an ample line bundle $L'$
on $X'$ such that $L \otimes([\f^{-1}(B)]) = \f^*L'$ and
$K_{X'} + (n-1)L'$ is ample.
\end{theorem}

A similar statement holds for the case $n=2$ where extra care is needed for
conic bundles.
Actually the adjunction theory and the above theorem in case of surfaces
were developed exactly a century ago
in a paper of Castelnuovo and Enriques (\cite{CaEn}, Firenze, fall 1900).\\
The above theorem was proved in the case $\ch(k) = 0$ by A. Sommese
who called the map $\f :X \to X'$ {\it the first reduction of $(X,L)$}
(see \cite{BS2} for a more general account as well as for applications of
adjunction theory).

We will actually prove a more general version of the theorem (see
(\ref{firstred})), where we consider
a rank $(n-1)$ ample vector bundle $\E$ and we study the non nef rays
in the cone $\overline{NE(X)}_{(K_X+det(\E)) \leq 0}$
(this version of the theorem if $\ch(k)= 0$ was proved in \cite{ABW2}).\par

\smallskip
In section 5 we restrict to the case $\ch(k) = 0$;
here we first generalize theorem (\ref{main}) with a
characterization of blow-up of smooth subvarieties in term of extremal rays
(see \ref{blow-up}).
In particular we prove a conjecture stated in \cite{AW}
for extremal rays whose associated contraction is divisorial;
we actually prove the conjecture for extremal rays which are non nef,
modulo normalization of the fiber (see (\ref{conj})).\\
Then we characterize non nef extremal rays of length $(n-2)$
(see (\ref{secred})). In particular this result generalizes the theorem which
guarantees the existence of the {\it second reduction of
a pair $(X,L)$} in the language of \cite[section 7.5]{BS2}.\par

\smallskip
The proof of our main theorem is strongly based on a recent result
by Y. Kachi and J. Koll\'ar \cite{Kako}; we report this result in
section 2.

\section{Set-up and preliminaries}

In this paper we will constantly use the notations and the definitions
presented in
the book of J\'anos Koll\'ar, \cite{Kob}, to which we will frequently refer.
For the reader's convenience we recall here the principal ones.

Let $X$ be a smooth $n$-dimensional projective variety
over an algebraically closed field of arbitrary characteristic.
A {\it rational curve} on $X$ is an irreducible reduced curve $C \subset X$
whose normalization is $f: \pu \to C \subset X$.

As we said in the introduction we assume that $K_X$ is not nef,
in particular that there exists
an {\it extremal ray} $R = \R^+[C]$ in the cone $\overline{NE(X)}_{{K_X} <0}$
(see \cite[III.1.2]{Kob}).
We will choose the curve $C$ generating the ray be to be rational
and minimal with respect to the intersection with $-K_X$;
sometimes we will say that such a $C$ is a {\it minimal extremal rational
curve}.
For such a curve $C$ the integer $l(R) = -{K_X}. C$ is called the {\it
length of the ray $R$}.

We will denote by $V \subset \om_{bir}(\pu,X)$ a closed irreducible
subvariety, which is closed under $\textrm{Aut}(\pu)$ and which contains $C$;
by the minimality of $C$ we have that $V$ is an {\it unsplit family of
rational curves} (see \cite[IV.2.]{Kob}).\\
We can define $\loc(V)$ to be the image in $X$
of the natural morphism associated to the family of rational curves $V$
(see \cite[II.2.3]{Kob}) and $\loc(V, 0 \to x)$ to be the image of
$V_x := \{f \in V: f(0) = x\}$; since $V$ is unsplit $\loc(V)$ is a proper
closed subscheme.\par
\medskip
We have the following fundamental inequality
\begin{proposition} \label{IW} \cite[IV.2.5 and 2.6.]{Kob} In the above
notation and for a general $x \in \loc(V)$
we have
$$\dim \loc(V) + \dim \loc(V, 0 \to x) +1 = \dim V \geq  \dim X +
deg_{(-{K_X})}(V).$$
\end{proposition}

We stress the fact that the equality holds since $V$ is an unsplit family.
The formulation and
the proof of this inequality is due to P. Ionescu and J. Wi\'sniewski, coming
out from the fundamental construction of S. Mori.

\smallskip
One can prove (see \cite[III.1.4]{Kob}) that, given an extremal ray $R$ on $X$,
there exists a Cartier divisor $H_R$ on $X$ such that
\begin{enumerate}
\item $H_R$ is nef.
\item If $z \in \overline{NE(X)}$ then ${H_R}. z = 0$ iff $z \in R$.
\end{enumerate}

Such an $H$ is called a {\it supporting divisor} for the extremal ray $R$.

\begin{proposition} \label{nefness}  (\cite[III.1.6]{Kob}). In the above
notation
one of the following two possibilities can occur
\begin{enumerate}
\item $(H_R)^{dimX} > 0$ iff there is an irreducible divisor $E \subset X$
such that $R.E < 0 $.
In particular $\loc(R)  \subset E$; in this case we say that the ray $R$
is {\sf not nef}.
\item $(H_R)^{dimX} = 0$ iff $Locus(R) = X$ (thus $X$ is uniruled).
In this case we say that the ray $R$ is {\sf nef}.
\end{enumerate}
\end{proposition}

Finally we recall the main result of the paper \cite{Kako} which we will
use in the proof of \ref{main}:
\begin{proposition} \label{kokacola} \cite[Lemma 8]{Kako} Let $X$ be a
normal, projective variety and $x\in X$
a smooth point. Let $H$ be an ample Cartier divisor on $X$ and let $\{C_t :
t\in T\}$
be an $(n-1)$-dimensional family of curves through $x$. Assume that
$(C_t.H)= 1$ and
$(-{K_X}.H^{n-1}) > (n-1) (H^n)$. Then $X \cong \pn$, $H = \opn (1)$ and
the $C_t$ are lines through $x$.
\end{proposition}

\section{Non nef extremal rays with maximal length}

In this section we prove the main theorem (\ref{main}) stated in the
introduction.

The if part of the theorem is trivial.\par
\smallskip
Assume therefore that $R$ is a non nef extremal ray of length $\geq (n-1)$;
thus by (\ref{nefness}) there exists an irreducible divisor
$E$ such that $\loc(R)$ is contained in $E$ and $E.R <0$.\\
Let $V$ be an unsplit family of deformations of $C$, a minimal extremal
rational curve generating the ray and let $x$ be a general point in $\loc(V)$;
applying the inequality in (\ref{IW}) we obtain that
$$\dim \loc(V) + \dim \loc(V, 0 \to x) +1 = \dim V \geq 2n -1.$$
On the other hand $\loc(V)\subseteq \loc(R) \subseteq E$ and thus
$$2n-2 \geq \dim \loc(V) + \dim \loc(V, 0 \to x),$$
forcing $\loc(V)=\loc(V, 0 \to x) = E$, $\dim V = 2n -1$ and $-{K_X}.C =
n-1$.\par
\smallskip
Let $f':\pu \to \Gamma \subset E$ be a curve of the family which intersects
the smooth
locus of
$E$; since every element of $\om_{[f']}(\pu,E)$ is also an element of
$\om_{[f']}(\pu,X)$ and $\loc (V) = E$ we can take an irreducible
component of $\om_{[f']}(\pu,E)$, call it $W$, which is contained in
$V$; this implies that $2n -1 = \dim V \ge \dim W$.\\
On the other hand, since $E$ is a locally complete intersection,
being a divisor in a smooth variety, we can apply \cite[Theorem II.1.3]{Kob}
which gives $\dim W \ge -K_E. \Gamma + n-1$; combining the two inequalities
we get
$$-K_E.\Gamma \le n.$$

\smallskip
Recalling that $E.R <0$, by the adjunction formula $K_E = (K_X +E)_{|E}$,
we have that $E.\Gamma = -1$ and $-K_E.\Gamma=n$.\par
\smallskip
Note that we can apply \cite[IV.3.13.3]{Kob} to the scheme $E$
and the family $V$; in fact for a generic point $x\in E$ we have seen that
$\loc(V, 0 \to x) = E$.\\
Thus we obtain that $\pic(E) = {\mathbb Z}$; moreover, from what we said above,
$-E_E$ is the positive generator of $\pic(E)$ and $-K_E = -nE_E$.\par
\smallskip
Let $\pi: \widetilde E \to E$ be the normalization of $E$; we have that
$$-K_{\widetilde E} = -\pi^*K_E +(\textrm{conductor of~} \pi)$$
where the conductor of $\pi$ is an effective divisor which is zero iff $\pi$
is an isomorphism (\cite[Ex. III.7.2]{Ha}).
Thus
$$-K_{\widetilde E }(\pi^*(-E_E))^{n-2} \geq {-K_E}(-E_E)^{n-2}
= n (-E_E)^{n-1}= n (\pi^*(-E_E))^{n-1}.$$
We lift the curves of $V$ to $\widetilde E$ and we apply \cite[Lemma 8]{Kako},
see \ref{kokacola},
which gives $\widetilde E = \proj^{n-1}$ and $\pi^*(-E_E) = \Ol_{\proj}(1)$.
Since
$$-K_{\widetilde E} = -\pi^*K_E +(\textrm{conductor of~} \pi) =
n \pi^*(-E_E) +(\textrm{conductor of~} \pi)$$ the conductor of $\pi$ is zero,
i.e. $\pi$ is an isomorphism.\par
\smallskip
To conclude the proof of the theorem we apply the following
\begin{theorem} (Castelnuovo) Let $E \subset X$ be a Cartier divisor in a
smooth projective variety $X$
with $E \cong \proj^{n-1}$ and $E_E = \Ol_{\proj}(-1)$; then there exists a
morphism
$\f : X \to X'$ into a smooth projective variety $X'$ which is the
blow-up of $X'$ at a point such that $E$ is the exceptional divisor.
\end{theorem}

The theorem was stated and proved in the case $\dim X = 2$ by Castelnuovo;
the same proof (see for instance \cite[V.5.7]{Ha}) applies in the general case.
\qed\par

\section{An application: the first reduction}

An immediate application of the main theorem gives the following
fundamental result for the so called {\it adjunction theory};
in the case $\ch(k) = 0$ the theorem was proved in \cite{ABW2}.

\begin{theorem}\label{firstred} Let $X$ be an $n$-dimensional smooth
variety over an algebraically closed field $ k$ with $n \geq 3$; let also
$\E$ be an ample vector bundle
on $X$ of rank $(n-1)$ such that $K_X +\det\E$ is nef and big, but not ample.\\
Then there exists a smooth variety $X'$ and a morphism
$\f :X \to X'$ expressing $X$ as blow-up of $X'$ at a finite set of points $B$
and an ample vector bundle $\E'$
on $X'$ such that $\E \otimes([\f^{-1}(B)]) = \f^*\E'$ and
$K_{X'} + \det\E'$ is ample.\\
The pair $(X',\E')$ is called the {\sf first reduction of $(X,\E)$}.
\end{theorem}

{\bf Proof.} \quad By Mori's cone theorem
there exists an extremal face $F \subset \overline{NE(X)}_{{K_X} <0}$
on which $K_X +\det \E$
is trivial; in other words $K_X +\det \E$ is the supporting divisor of a face
$F$.
By assumption every extremal ray in $F$ is non nef
and has length $\geq (n-1)$ (use propositions \ref{nefness} and \ref{IW}),
hence we can apply the above theorem
to each of these rays.\\
If $R_i$ and $R_j$ are two different rays in $F$ then the respective
Loci, $E_i$ and $E_j$, are disjoint. In fact, assume the contrary;
then, by Serre's inequality, we have

$$\dim (E_i \cap E_j) \ge \dim E_i +\dim E_j -n \ge n-2 \ge 1$$

so that the intersection of $E_i$ and $E_j$ would contain a curve $B$
and this is a contradiction.\\

Finally let $\E'$ be a rank $n-1$ vector bundle such that
$\E \otimes(\sum E_i) = \f^*\E'$; it is easy to check that $\E'$ and $K_{X'}
+\det \E'$
are ample, using the Kleiman ampleness criterium, and that $K_X +\det \E =
\f^*(K_{X'}+\det \E')$. \qed \par
\medskip
This result implies the one stated in the introduction taking $\E =
\oplus^{(n-1)} L$
and using the following simple result:

\begin{lemma} Let $X$ be an $n$-dimensional smooth
variety over an algebraically closed field $ k$ with $n \geq 3$; let also
$\E$ be an ample vector bundle
on $X$ of rank $\geq (n-1)$.
If there is a non nef ray in the cone
$\overline{NE(X)}_{(K_X+ det \E) \leq 0}$ then,
with the only exception given by the blow-up of $\pt$
in one point, $\pi : Bl_x \pt \to\pt$, and $\E = \oplus^2 (\pi^*(\opt (2))
-[\pi^{-1}(x)])$,
all rays in this cone are non nef.
\end{lemma}

{\bf Proof.} \quad Assume we have a non nef ray $R_1$ and a nef one $R_2$ as in
the lemma.
The ray $R_1$ satisfies the assumptions of theorem (\ref{main}), thus in
particular $\loc
(R_1) =\pnm$.
Let $V_2 \subset \om_{bir}(\pu,X)$ be the family associated to $R_2$
described in section 2 and let $x \in \loc(R_1)$ be a general point; since by
(\ref{IW}) $\dim \loc(V_2, 0 \to x) \geq n-2$  we have, by Serre's
inequality, that
the intersection $\loc(V_2, 0 \to x) \cap \loc(R_1)$ contains a curve, unless
$n=3$ and $\dim \loc(V_2, 0 \to x) =1$.\\
In the general case this is a contradiction;
in the special one, since extremal rays and their contractions are known in
all characteristic
in dimension $3$ (see \cite{Ko0}), we have that the contraction of the
extremal ray
$R_2$, $\f:X \to Z$ gives $X$ a structure of a $\pu$ bundle
of which $\loc(R_1)=\pd$ is a (multi)section.
Let $\pi:X \to X'$ be the contraction of $R_1$ and $x' = \pi(\loc(R_1))$;
the images of the fibers of $\f$ via $\pi$ are irreducible rational curves
passing through $x'$ and covering $X'$, hence we can apply
\cite[IV.3.13.3]{Kob}
to deduce $\pic(X')= \mathbb Z$.
By the blow up formula $-K_{X'}$ is ample so $X'$
is Fano; then one can compute that the length of the unique
ray of $X'$ is $4$ and therefore that $X' \simeq \pt$, by the classification
of Fano threefolds \cite{SB}.\qed\par

\begin{remark}\label{nef} The description of the nef rays in
$\overline{NE(X)}_{(K_X+ (n-1)L) \leq 0}$
(or, more generally, in
$\overline{NE(X)}_{(K_X+ det \E) \leq 0}$)
is well understood in characteristic zero (\cite{Fu1} and \cite{ABW2}).
In positive characteristic the description of such rays was recently given if
$l(R) > (n-1)$, see \cite{Kako}.\\
If $l(R) = n-1$ we cannot at the moment say too much,
the main problem being the lack of a general contraction theorem
for Mori rays.\\
Studying the families of rational curves
arising from extremal rays we can prove that, if there exists
a ray of this kind, either $\pic(X) \simeq \mathbb Z$ and
$-K_X = (n-1)H$ or
there exists an unsplit covering family $V$ of rational curves
such that $\dim \loc(V, 0 \to x) = n-2$, for a general $x \in X$.
\end{remark}

\section{Non nef extremal rays of high length}

In this section we restrict ourselves to characteristic zero.
In this case to each extremal ray $R$ is associated an extremal
contraction $\f_R : X \to X'$, that is a morphism onto a normal projective
variety
$X'$ with connected fibers and such that if $C \subset X$ is a curve then
$\f(C) $ is a point if and only if $[C] \in R$.
Denoting by $E(\f)$ the exceptional locus of $\f$ and by $S$ an
irreducible component of a (non trivial) fiber, the inequality in
\ref{IW} can be read as
$$\dim E(\f) + \dim S \geq \dim X + l(R) -1.$$

We first prove that the same argument in the proof of theorem
(\ref{main}) can be used to (almost) prove Conjecture (2.6) in \cite{AW}.
Namely we have the following two theorems:

\begin{theorem}\label{conj} Let $X$ be a smooth $n$-dimensional
projective variety
over an algebraically closed field of characteristic zero
and $R=\mathbb R_+[C]$ a non nef extremal ray.
If the above inequality holds for an irreducible component $S$ then
its normalization, $\rho : S' \to S$, is ${\mathbb P}^s$.
\end{theorem}

In the following more special case we have the optimal characterization:

\begin{theorem}\label{blow-up} Let $X$ be a smooth $n$-dimensional
projective variety
over an algebraically closed field of characteristic zero; the two following
facts are equi\-va\-lent:\\
1) There exists an extremal ray $R$ such that the contraction associated to
$R$ is divisorial
and the fibers have dimension $= l(R)$. \\
2) There exists a morphism
$\f : X \to X'$ into a smooth projective variety $X'$ which is the
blow-up of $X'$ along a smooth subvariety of codimension $l(R) +1$.
\end{theorem}

{\bf Proof.} \quad We will only prove that 1) implies 2) in the theorem
\ref {blow-up},
the other implication being trivial.
We observe also that for both theorems
it is enough to
prove that there exists a line bundle $L$ on $X$ such that $L.C=1$;
in fact, if this is the case, the contraction of $R$ will be supported by a
divisor
of the type $K_X+l(R) L$ with $L$ a $\f$-ample line bundle. Then, theorem
\ref{conj}
will follows from \cite[Lemma 1.1]{ABW1} while theorem \ref{blow-up}
will follow from \cite[Theorem 4.1.iii]{AWD}.

Let then $R$ be an extremal ray as in \ref{conj} or as in 1) of \ref{blow-up}.
Since $R$ is not nef, there exists an irreducible divisor
$E$ such that $E(\f) = \loc(R)$ is contained in $E$ and $E.R <0$;
let $V$ be an unsplit family of deformations of $C$.
Since by the assumptions of both theorems the inequality in (\ref{IW}) is
an equality
we have that
$$\dim V = n + l(R) = n  -{K_X}.C.$$
On the other hand we can think at $V$ as a family in $E$, as in the proof of
\ref{main}, thus, if we apply \cite[Theorem II.1.3]{Kob},
we have $\dim V \ge -K_E.C + n-1$.
Combining the two inequalities
we get
$$-K_E.C \le l(R)+1.$$

Recalling that $E.R <0$ and that by the adjunction formula
$K_E = (K_X +E)_{|E}$,
we have that $E.C = -1$.\qed\par

\medskip
In the special case of characteristic zero we can further extend
our classification of non nef extremal rays with high
(but not maximal) length as in the following theorem.

\begin{theorem} \label{secred} Let $X$ be a smooth $n$-dimensional
projective variety
over an algebraically closed field of characteristic zero
and $R=\mathbb R_+[C]$ a non nef extremal ray.
Let $\f:X \to X'$ be the birational elementary contraction associated to $R$
and $E=\textrm{Exc}(\f)$ its exceptional locus.
If $l(R) = n-2$ then one of the following cases occur:\\
1)  $\f(E)$ is a point and $(E,-E_E) \simeq (\proj^{n-1}, \Ol_{\proj}(2))$.\\
2)  $\f(E)$ is a point and $(E,-E_E) \simeq
(\mathbb Q^{n-1}, \Ol_{\mathbb Q}(1))$, where $\mathbb Q^{n-1}$ is a
possibly singular quadric.\\
3) $X'$ is smooth and $\f$ is the blow-up along a smooth curve $\f(E)
\subset X'$.
\end{theorem}

{\bf Proof.} \quad Also here it is enough to
prove that there exists a line bundle $L$ on $X$ such that $L.C=1$;
then the theorem will follow from \cite[Theorem 4]{Fu1}
and \cite[Theorem 4.1.iii]{AWD}.

Since $R$ is not nef, there exists an irreducible divisor
$E$ such that $\loc(R)$ is contained in $E$ and $E.R <0$.\\
Let $V$ be the unsplit family of deformations of $C$, a minimal extremal
rational
curve; by the inequality in (\ref{IW}) we have
that $\loc(R)=E$  and  $\dim V =\dim \loc(V) + \dim \loc(V, 0 \to x) +1$
is either $2n-1$ or $2n-2$.\\
On the other hand we can think at $V$ as a family in $E$, as in the previous
proofs. Thus, if we apply \cite[Theorem II.1.3]{Kob},
we have $\dim V \ge -K_E.C + n-1$.
In particular we get
$$-K_E.C \le n.$$

Recalling that $E.R <0$ and that by the adjunction formula
$K_E = (K_X +E)_{|E}$,
we have that $E.C = -1$ or $-2$.

In the first case we are done;
the second one, on the other hand, can occur only if $\dim V =2n-2$
and $\dim \loc(V, 0 \to x) = n-1$.
In particular we have that $\loc(V, 0 \to x) = E$, so we can apply
\cite[IV.3.13.3]{Kob} to get $\pic(E)\simeq \mathbb Z$.
Let $\pi:\widetilde E \to E$ be the normalization of $E$;
lifting the family $V$ to $\widetilde E$ we obtain a new family which,
together with $-\pi^*E_E$,
satisfies on $\widetilde E$ the assumptions of \cite[Theorem 3.6]{Keb}.
This yields that
$\widetilde E \simeq \proj^{n-1}$ and $-\pi^*E_E \simeq \Ol_{\proj}(2)$.
Let $H_E$ be the ample generator of $\pic(E)$; we can write
$-E_E = e H_E, \quad -K_E= ({{ne} \over 2})H_E, \quad \pi^*H_E = \Ol_{\proj}(2/e)$; then
$$-K_{\widetilde E} = -\pi^*K_E + (\textrm{conductor of $\pi$}),$$
hence 
$${{ne}\over 2}\Ol_{\proj}(2/e) = -K_{\widetilde E}= \pi^*({{ne} \over 2})H_E +
(\textrm{conductor of $\pi$})$$
so that the conductor is zero and $\pi$ is an isomorphism.
\qed\par

\par
\medskip
{\small Marco Andreatta \\
Dipartimento di Matematica \\
Universit\`a degli Studi di Trento \\
Via Sommarive 14\\
I-38100 Povo (TN), Italy\\
e-mail: andreatt@science.unitn.it\\
Fax: +39-0461-881624}\par
\medskip
{\small Gianluca Occhetta \\
Dipartimento di Matematica ``F. Enriques''\\
Universit\`a degli Studi di Milano \\
Via Saldini, 50\\
I-20133 Milano, Italy\\
e-mail: occhetta@mat.unimi.it}

\end{document}